\documentstyle[12pt]{article}
\def\lb{\label}
\def\be{\begin{equation}}
\def\ee{\end{equation}}

\begin{document}

\begin{center}
{\large \bf BRST OPERATOR FOR QUANTUM LIE ALGEBRAS:
RELATION TO BAR COMPLEX}

{\it V. Gorbounov$^*$, A.P. Isaev$^{**}$ and O. Ogievetsky$^{***}$}

{\it $^*$Department of Mathematics, University of Kentucky, USA \\
email: vgorb@ms.uky.edu} \\
{\it $^{**}$ Bogoliubov Laboratory of Theoretical
Physics, Joint Institute for Nuclear Research, Dubna 141980,
Moscow reg., Russia \\
email: isaevap@thsun1.jinr.ru} \\
{\it $^{***}$ Center of Theoretical Physics, Luminy,
13288 Marseille, France
and P. N. Lebedev Physical Institute, Theoretical Department,
Leninsky pr. 53, 117924 Moscow, Russia \\
email: oleg@cpt.univ-mrs.fr}

\end{center}

\begin{abstract}
Quantum Lie algebras
(an important class of
quadratic algebras arising in the Woronowicz calculus on quantum groups)
are generalizations of Lie (super) algebras.
Many notions from the theory of Lie (super)algebras admit
``quantum'' generalizations.
In particular, there is a BRST operator $Q$
($Q^2=0$) which generates the differential in the Woronowicz theory
and gives information about (co)homologies of quantum Lie algebras.
In our previous papers a recurrence relation for the operator $Q$
for quantum Lie algebras was given and solved.
Here we consider the bar complex for q-Lie algebras and
its subcomplex of q-antisymmetric chains.
We establish a chain map (which is an isomorphism) of the standard complex
for a q-Lie algebra to
the subcomplex of the antisymmetric chains.
The construction requires a set of nontrivial identities in the group algebra
of the braid group. We discuss also a generalization of the
standard complex to the case
when a q-Lie algebra is equipped with a grading operator.
\end{abstract}




\section{Introduction}

\vspace{0.3cm}

The Woronowicz calculus \cite{Wor} associates
an algebra of exterior forms $\Gamma$ and an enveloping algebra $\mathcal U$
of the left invariant vector fields
on $\mathcal A$  to a Hopf algebra $\mathcal A$. 
The algebra ${\mathcal U}$ is called a quantum Lie algebra
(q-Lie algebra for short). It is defined by relations
$\chi_{i} \, \chi_{j} - \sigma^{km}_{ij} \, \chi_{k} \, \chi_{m} =
C_{ij}^{k} \, \chi_k$, where $\{ \chi_i \}$ is a set of generators;
"structure constants" $\sigma^{km}_{ij}$ and $C_{ij}^{k}$ obey certain
constraints, see Section 2. A q-Lie algebra is a non-homogeneous quadratic algebra.
The general theory of quadratic algebras has been considered in
a number of papers (see e.g. \cite{Zyg,BG,PP}). The case of q-Lie algebras
is quite particular. An analog of the de Rham complex
for $\mathcal U$ has been constructed in \cite{Wor}.
Mimicking the classical theory of Lie algebras,
quantum analogs of the standard complex and BRST differential (with expected properties)
have been introduced in \cite{isog,big}
(for a review of the general BRST-theory see \cite{HT}).
In this paper we continue the investigation of
q-Lie algebras along the lines of the theory of Lie algebras.
Namely, (see e.g. \cite{CE}) we map the standard complex of the q-Lie algebra into
the subspace of q-antisymmetric chains
of the bar complex for $\mathcal U$. It turns out that the image of this map,
just like in the case Lie algebras, is a subcomplex
in the bar complex for $\mathcal U$. Moreover, the
quantum BRST differential is the restriction of the bar differential $b$ to
the subspace of the q-antisymmetric chains.

The paper is organized as follows. In Section 2 we recall
the definition of q-Lie algebras.
In Section 3 we introduce the exterior
extensions of q-Lie algebra (an algebra of exterior forms $\Gamma$ and
inner derivatives $\Gamma^*$)
due to Woronowicz (our definition of the
basis of differential
forms is slightly different from the Woronowicz definition, see \cite{isog}) and develop the
tools needed for the main constructions of Sections 4 and 5.
Section 4 contains an explicit
construction of the BRST operator of a q-Lie algebra. The presentation in
the Sections 2-4 follows \cite{isog2}.

In Section 5 we discuss the bar complex  $(C_{n}({\mathcal U}),b)$
for a q-Lie algebra,
its subcomplex of q-antisymmetric chains $(C_{n}({\mathcal U}, \wedge {\mathcal L}),b)$
and a chain map of the standard complex for q-Lie algebra to
$(C_{n}({\mathcal U}, \wedge {\mathcal L}),b)$.
Nontrivial identities in braid group algebra are
used for proving the main Propositions of Sections 4 and 5.

Section 6 is a generalization of the constructions above to the situation
when a q-Lie algebra is equipped with a grading operator. An example:
a Lie super-algebra with grading given by parity. It is known that there
are two choices of commutation relations between bosonic and fermionic ghosts.
We explain this phenomenon in a general framework of a q-Lie algebra
with a grading operator.

Section 7 summarizes the results of the paper.

\section{Definitions and notation}

The data defining a quantum Lie algebra with $N$ generators can be
conveniently encoded
into the following $(N+1)^2 \times (N+1)^2$ matrix \cite{Ber}:
\be
\label{rmat}
\begin{array}{c}
R^{CD}_{AB}
 = \delta^{i}_A \, \delta^{j}_B \, \sigma_{ij}^{kl} \,
\delta_{k}^C \, \delta_{l}^D +
\delta^{i}_A \, \delta^{j}_B \, C_{ij}^{k} \,
\delta_{0}^C \, \delta_{k}^D + \delta_A^{i} \delta^D_{i} \,
\delta^{0}_B \, \delta_{0}^C
+ \delta^C_B \, \delta^{0}_A \, \delta_{0}^D   \\ \\
\equiv \delta^{\langle 1}_A \, \delta^{\langle 2}_B \, \sigma_{12} \,
\delta_{1\rangle }^C \, \delta_{2\rangle }^D +
\delta^{\langle 1}_A \, \delta^{\langle 2}_B \, C_{12\rangle }^{\langle 2} \,
\delta_{0}^C \, \delta_{2\rangle }^D + \delta_A^{\langle 1}
\delta^D_{1\rangle } \, \delta^{0}_B \, \delta_{0}^C
+ \delta^C_B \, \delta^{0}_A \, \delta_{0}^D   \; .
\end{array}
\ee
Here $1,2, \dots$ denote copies of an $N$-dimensional vector spaces $V_N$;
$\langle  1 , \langle  2 , \dots$
(resp., $1 \rangle  , 2 \rangle  , \dots$) are the corresponding outcoming
and incoming vectors; small letters
$i,j,k,l,\dots = 1,2, \dots , N$ denote indices of the vector space
$V_N$; capital letters
 $A,B, \dots = 0, 1, \dots , N$ denote indices of an
$(N+1)$-dimensional space $V_{N+1}$
($V_N$ is a subspace in $V_{N+1}$).

Eqs. (\ref{rmat}) are equivalent to
\be
\label{rmat1}
R^{ij}_{kl} = \sigma^{ij}_{kl} \; , \;\;\; R^{0j}_{kl} = C^{j}_{kl} \; , \;\;\;
R^{0A}_{B0} = R^{A0}_{0B} = \delta^A_B \;
\ee
and the other components of $R$ vanish. We assume that the
matrix $\sigma\in End(V_{N} \otimes V_{N})$ has an eigenvalue $1$.

Let $L_A = \{ \chi_0,\chi_i \} $ $(i = 1, \dots, N)$ be a
``quantum vector'':
\be
\label{qv}
L_A =  \delta^{\langle 1}_A \chi_{1\rangle }  + \delta^0_A \, \chi_0  =
\left(
\begin{array}{c}
\chi_0 \\
\chi_i
\end{array}
\right) \;
\ee
for the matrix $R$:
\be
\label{rtt1}
R^{CD}_{AB} \, L_C \, L_D = L_A \, L_B \; .
\ee
These relations are equivalent to
$$
[\chi_0 , \, \chi_i] = 0\ \ {\rm and}\ \
(1 - \sigma_{12}) \, \chi_{1\rangle } \, \chi_{2\rangle } =
C_{12\rangle }^{\langle 1} \, \chi_0 \, \chi_{1\rangle }\ .
$$
One can rescale
the elements, $\chi_i \rightarrow \chi_0 \, \chi_i$. The algebra with
rescaled generators
$\chi_i$ ($i =1,2, \dots, N$) subject to relations
\be
\label{qla}
(1 - \sigma_{12}) \, \chi_{1\rangle } \, \chi_{2\rangle } =
C_{12\rangle }^{\langle 1} \, \chi_{1\rangle } \; ,
\ee
is called a {\it quantum Lie algebra} if
the matrix $R$ satisfies the Yang-Baxter equation,
\be
\label{rrr}
R^{CD}_{AB} \, R^{ES}_{MC} \, R^{FN}_{SD} = R^{CS}_{MA} \, R^{DN}_{SB} \,
R^{EF}_{CD}
\ \ {\rm or}\ \
R_{\underline{2} \underline{3}} \, R_{\underline{1} \underline{2}} \,
R_{\underline{2} \underline{3}} =
R_{\underline{1} \underline{2}} \, R_{\underline{2} \underline{3}} \,
R_{\underline{1} \underline{2}} \;
\ee
(here $\underline{1},\underline{2}$ or $\underline{2},\underline{3}$ denote copies of the
vector spaces $V_{N+1}$ on which $R$- matrices
act nontrivially). The Yang-Baxter equation for $R$ imposes the
following conditions for the structure constants of the quantum Lie algebra
(\ref{qla}):
\be
\label{rell1}
\sigma_{1} \, \sigma_{2} \, \sigma_{1} =
\sigma_{2} \, \sigma_{1} \, \sigma_{2}
\; ,
\ee
\be
\label{rell2}
C_{1} \, \delta_3 \, C_{1} =
\sigma_{2} \, C_{1} \, \delta_3 \, C_{1} +
C_{2} \, C_{1} \; ,
\ee
\be
\label{rell3}
\vspace{2mm}
C_{1} \, \delta_3 \, \sigma_{1} = \sigma_{2} \, \sigma_{1} \, C_{2} \; ,
\ee
\be
\label{rell4}
\begin{array}{c}
\vspace{2mm}
(\sigma_{2} \, C_{1} \, \delta_3 +C_{2}) \, \sigma_{1}=
\sigma_{1} \, (\sigma_{2} \, C_{1} \, \delta_3 + C_{2})\ .
\end{array}
\ee
Here we have used a concise notation $\sigma_n := \sigma_{n \, n+1}$,
$C_n := C^{\langle n}_{n,n+1\rangle }$, $\delta_n := \delta_{\;\;
n\rangle }^{\langle n-1}$.

\section{Exterior extensions for quantum Lie algebra}

To define the exterior extension of the quantum Lie algebra (\ref{qla}),
one introduces quantum wedge algebras with
generators $\gamma_i$ and $\Omega^i$ ($i =1, \dots N$):
 \be
 \label{owedge}
\Omega^{\langle n}  \wedge \dots \wedge \Omega^{\langle 2} \wedge
\Omega^{\langle 1}
= \Omega^{\langle n} \otimes  \dots \otimes
\Omega^{\langle 2} \otimes \Omega^{\langle 1} \, A_{1 \rightarrow n} \; ,
 \ee
  \be
 \label{gwedge}
\gamma_{1\rangle } \wedge \gamma_{2\rangle } \dots \wedge \gamma_{n\rangle } =
A_{1 \rightarrow n} \, \gamma_{1\rangle } \otimes \gamma_{2\rangle } \dots
\otimes \gamma_{n\rangle } \; .
 \ee
The cross-commutation relations are
\be
\label{go}
 \gamma_{2\rangle } \, \Omega^{\langle 2} = - \Omega^{\langle 1} \,
\sigma^{-1}_{12} \, \gamma_{1\rangle } + I_2  \; .
\ee
In eqs. (\ref{owedge}) and (\ref{gwedge}) we have used operators
$A_{1 \rightarrow n}$, quantum analogues
of the antisymmetrizers. They can be defined inductively
(in fact they make sense already for the group algebra of the braid
group $B_{M+1}$ with generators $\sigma_i$ $(i=1, \dots, M)$):
\be
\label{anti}
A_{1\to n} \equiv  f_{1\to n} \, A_{1\to n-1} \ \ {\rm or}\ \
A_{1\to n} \equiv  \overline{f}_{1\to n} \, A_{2\to n}   \; ,
\ee
where $A_{1\to 1}  = 1$,  $f_{k \to k} = \overline{f}_{k \to k}=1$ and
\be
\label{fk}
\begin{array}{l}
f_{k \to m} = 1 - f_{k \to m-1} \, \sigma_{m-1} \ \ {\rm or}\ \
f_{k \to m} = f_{k +1 \to m} + (-1)^{m-k} \, \sigma_{k} \cdots \sigma_{m-1}
 \; ,
\\ \\
\overline{f}_{k \to m} = 1- \overline{f}_{k+1 \to m} \, \sigma_{k}
\ \ {\rm or}\ \
\overline{f}_{k \to m} = \overline{f}_{k \to m-1}
+ (-1)^{m-k} \, \sigma_{m-1} \cdots \sigma_{k+ 1} \, \sigma_{k}  \;
\end{array}
\ee
for $(k < m \leq n)$. Explicitly:
\be
\label{fbf}
\begin{array}{l}
f_{k \to n} \equiv
f^{(\sigma)}_{k \to n} =
 1 - \sigma_{n-1} + \sigma_{n-2} \,
\sigma_{n-1} - \dots
+ (-1)^{n-k} \, \sigma_{k} \, \sigma_{k+ 1} \cdots \sigma_{n-1}  \;  , \\ \\
\overline{f}_{k \to n} \equiv
\overline{f}^{(\sigma)}_{k \to n} =
1 - \sigma_{k} + \sigma_{k+1} \, \sigma_{k} - \dots
+ (-1)^{n-k} \, \sigma_{n-1} \cdots \sigma_{k+ 1} \, \sigma_{k}\; .
\end{array}
\ee
If the sequence of operators $A_{1 \rightarrow n}$ terminates
at the step $n = h+1$ ($A_{1 \rightarrow n} = 0$ $\forall n >h$)
for some matrix representation $\rho$ of $B_{M+1}$
then the number $h$ is called the height of
the matrix $\rho(\sigma_{12})$.

The elements $\gamma_i$ and $\Omega^i$ are quantum analogues of
ghost variables.
The choice of the sign and the appearance of the braid matrix
in (\ref{go}) follows the conventional case of Lie (super)algebras.
The commutation relations of $\chi_i$ with $\gamma_j$ and
$\Omega^k$ are
\be
\label{cross}
\chi_{2\rangle } \, \Omega^{\langle 2} =  \Omega^{\langle 1} \, \left(
\sigma_{12} \, \chi_{1\rangle } +  C^{\langle 2}_{12\rangle } \right) \, ,
\;\; \gamma_{1\rangle } \, \chi_{2\rangle } = \sigma_{12} \, \chi_{1\rangle } \,
\gamma_{2\rangle } + C^{\langle 1}_{12\rangle } \, \gamma_{1\rangle } \; .
\ee

The algebra (\ref{qla}), (\ref{owedge}) -- (\ref{go}) and (\ref{cross})
is the exterior extension of the quantum Lie algebra and
(see \cite{isog}) this
algebra gives rise to the Cartan differential calculus on quantum groups
in the Woronowicz theory \cite{Wor}.

\vspace{0.5cm}

\noindent
{\bf Remark.}
Let $\hat{\sigma}_i$ $(i=1, \dots M)$ be generators of
the braid group $B_{M+1}$:
\be
\label{braidg}
\hat{\sigma}_i \, \hat{\sigma}_{i+1} \, \hat{\sigma}_i =
\hat{\sigma}_{i+1} \, \hat{\sigma}_i \,  \hat{\sigma}_{i+1} \; , \;\;\;
[ \hat{\sigma}_{i},  \hat{\sigma}_{j}] = 0 \;
{\rm for} \; |i-j| > 1 \; .
\ee
The important properties of the elements $f^{(\hat{\sigma})}_{k \to n}$ and
$\overline{f}^{(\hat{\sigma})}_{k \to n}$ (\ref{fbf}) of $B_{M+1}$ are:
\be
\lb{ffxa}
f^{(\hat{\sigma})}_{1 \to n} f^{(\hat{\sigma})}_{1 \to n-1}
\dots f^{(\hat{\sigma})}_{1 \to m} = x^{(n-m+1)}_{n} \, A_{m \rightarrow n} \; , \;\;\;
(n \geq m \geq 1)
\ee
\be
\lb{ffxb}
\overline{f}^{(\hat{\sigma})}_{1 \to n} \overline{f}^{(\hat{\sigma})}_{2 \to n}
\dots \overline{f}^{(\hat{\sigma})}_{m \to n} =
x^{(n-m)}_{n} \, A_{1 \rightarrow m} \; , \;\;\;
(n \geq m \geq 1)
\ee
where $x^{(0)}_{n} =1$, $x^{(1)}_{m} = f^{(\hat{\sigma})}_{1 \to m}$,
$x^{(n)}_{n} =1$, $x^{(n-1)}_{n} =\overline{f}^{(\hat{\sigma})}_{1 \to n}$
and $x^{(2)}_{m + 1}$ is given by
\be
\lb{xm}
\begin{array}{c}
x^{(2)}_{m+1} = f^{(\hat{\sigma})}_{1 \to m} +
f^{(\hat{\sigma})}_{1 \to m-1} \,
(\hat{\sigma}_{m} \, \hat{\sigma}_{m-1}) +
f^{(\hat{\sigma})}_{1 \to m-2} \, (\hat{\sigma}_{m-1} \,
\hat{\sigma}_{m-2}) \,
(\hat{\sigma}_{m} \, \hat{\sigma}_{m-1}) \\ \\
+ \dots + f^{(\hat{\sigma})}_{1 \to 2} \, (\hat{\sigma}_{3} \,
\hat{\sigma}_{2}) \cdots
(\hat{\sigma}_{m} \, \hat{\sigma}_{m-1})
+ (\hat{\sigma}_{2} \, \hat{\sigma}_{1}) \cdots (\hat{\sigma}_{m} \,
\hat{\sigma}_{m-1}) \; .
\end{array}
\ee
The identities (\ref{ffxa}), (\ref{ffxb}) are equivalent to the factorization formula
\be
\lb{shuffle}
A_{1 \to n} = x^{(n-m)}_{n} \, A_{1 \to m} \, A_{m +1 \rightarrow n}
\ee

The elements $x^{(m)}_{n}$ can be defined inductively using the recurrent relation
$$
 x^{(n-m+1)}_{n+1} = x^{(n-m)}_{n}  - (-1)^{n-m} \,
 x^{(n-m+1)}_{n} \, \hat{\sigma}_{n} \dots \hat{\sigma}_{m}
$$
and in particular we have
\be
\lb{indux}
x^{(2)}_{n+1} =  (f^{(\hat{\sigma})}_{1 \to n} + x^{(2)}_{n} \,
\hat{\sigma}_{n} \, \hat{\sigma}_{n-1}) \; , \;\;\;
(n \geq 2) \; ,
\ee
which gives (\ref{xm}). Then, we have
$x^{(3)}_{n+1} =  x^{(2)}_{n}  - x^{(3)}_{n} \, \hat{\sigma}_{n} \,
\hat{\sigma}_{n-1} \, \hat{\sigma}_{n-2}$
for $(n \geq 3)$, etc.

Note that the elements $x^{(m)}_{n}$ are alternating sums over the braid
group elements which can be considered as quantum analogs of $(m, \, n-m)$
shuffles (the shuffles are obtained by projection $\hat{\sigma}_i \to
-s_i$, where $s_i$ are generators of symmetric group $S_{M+1}$).
From this point of view the wedge products (\ref{owedge}), (\ref{gwedge})
are related to the quantum shuffle products
(about quantum shuffles and corresponding products see \cite{Rosso},
\cite{Wor}). The associativity of these
products is provided by the identities: $x_n^{(n-m)} \, x_m^{(m-k)} =
x_n^{(n-k)} \, (T_k \, x_{n-k}^{(n-m)} \, T_k^{-1})$ for $(k < m < n)$,
where $T_k \, \hat{\sigma}_i \, T_k^{-1} = \hat{\sigma}_{i+k}$.


\section{BRST operator for quantum Lie algebra}

In the paper \cite{isog} we have given a recurrence which determines
the BRST operator $Q$ for the quantum Lie algebra:

\vspace{0.3cm}
\noindent
{\bf Proposition 1.}
{\it The BRST operator $Q$ for the quantum Lie algebras (\ref{qla}),
which satisfies the equation $Q^2=0$, should have the following form
\begin{equation}
\label{brst}
Q = \Omega^i \, \chi_i + \sum^{h-1}_{r=1} \, Q_{(r)} \; ,
\end{equation}
where $h$ is the height of the braid matrix $\sigma_{12}$,
the operators $Q_{(r)}$ are given by
\begin{equation}
\label{qrr}
Q_{(r)} =\Omega^{\langle r+1|} \otimes \Omega^{\langle r|} \otimes \dots
\otimes \Omega^{\langle 1|} \, (AXA)^{\;\; \langle 1 \dots r|}_{|1 \dots
r+1\rangle } \,
\gamma_{| 1 \rangle } \otimes \dots \otimes \gamma_{| r \rangle }
\end{equation}
and $(AXA)^{\;\; \langle 1 \dots r |}_{|1 \dots r+1\rangle }
= A_{1 \to r+1}( XA)^{\;\; \langle 1 \dots r |}_{|1 \dots r+1\rangle } =
(AX)^{\;\; \langle 1 \dots r |}_{|1 \dots r+1\rangle} A_{1 \to r}$ are tensors
which obey the following recurrent
relations
\begin{equation}
\label{recur}
A_{1 \rightarrow r+1} \,
(XA)^{\;\; \langle  1\dots r|}_{|1 \dots r+1\rangle }  =
A_{1 \rightarrow r+1} \, \left((-1)^{r} \,
\sigma_{r} \sigma_{r-1} \cdots  \sigma_{1} - {\bf 1} \right) \,
(XA)^{\;\; \langle  2\dots r |}_{|2 \dots r+1\rangle }  \;
\end{equation}
with the initial condition $A_{12} \, X^{\langle 1|}_{|12\rangle } = -
C^{\langle 1|}_{|12\rangle }$.
  }

\vspace{.3cm}

It should be stressed that the derivation of the recurrence (\ref{recur})
requires the consideration of the terms (in the equation $Q^2=0$) linear in $\chi_i$ only.

\vspace{.3cm}

In \cite{isog2} we presented the solution of eqn. (\ref{recur}) and defined
coefficients \\
$(A \, X \, A)^{\;\; \langle  1\dots r|}_{|1 \dots r+1\rangle }$ explicitly.
This was done with the help of an important identity in the
group algebra of the braid group:

\vspace{0.3cm}
\noindent
{\bf Lemma 1.} {\it Let $\hat{\sigma}_i$ $(i=1, \dots M)$ be generators of
the braid group $B_{M+1}$ (\ref{braidg}).
Then the following identity in the group algebra of $B_{M+1}$ holds:
\be
\label{rel15}
Y_{1 \to r+1}
 \, \overline{f}^{(\hat{\sigma})}_{1 \to r}
=  \left((-1)^{r+1} \, f^{(\hat{\sigma})}_{1\to r+1 } \,
 \hat{\sigma}_{r}  \cdots  \hat{\sigma}_{1} +
 \overline{f}^{(\hat{\sigma})}_{1\to r+1} \right)
\, Y_{2 \to r+1} \; ,
\ee
where $f^{(\hat{\sigma})}_{1\to r}$ and $\overline{f}^{(\hat{\sigma})}_{1 \to r}$
are defined in (\ref{fbf}) and
\be
\label{yy}
\begin{array}{l}
Y_{k \to r+1} =
(1-\hat{\sigma}^2_{r}) \left( 1 + \hat{\sigma}_{r-1} \,
\hat{\sigma}^2_{r} \right)
\cdots \left( 1 - (-)^{r-k} \, \hat{\sigma}_{k} \,
\hat{\sigma}_{k+1} \dots \hat{\sigma}_{r-1} \, \hat{\sigma}_{r}^2 \right) \; ,
\;\;\; \\[1em]
Y_{k \to k} = 1 \; .
\end{array}
\ee
}
\vspace{0.1cm}
\noindent
{\bf Proof.} With the help of eqs. (\ref{fk}), 
the identity (\ref{rel15}) can be represented in the form
\be
\label{rel16}
Y_{1 \to r+1}
 \, \overline{f}^{(\hat{\sigma})}_{1 \to r}
= F_{1 \to r+1} \, Y_{2 \to r+1}\ ,
\ee
where
\be
\label{Ff}
F_{1 \to r+1} = \left((-1)^{r} \, f^{(\hat{\sigma})}_{1\to r } \,
 \hat{\sigma}^2_{r}  \cdots  \hat{\sigma}_{1} +
 \overline{f}^{(\hat{\sigma})}_{1\to r} \right) \; .
\ee
One can prove the identity (\ref{rel16}) by induction
(see \cite{isog2}, \cite{MG}). $\bullet$

\vspace{0.3cm}
\noindent
{\bf Remark.}
The operators
$Y_{k \to r+1}$ (\ref{yy}) have another form:
\be
\label{zag}
Y_{k \to r+1} =
(1 + \hat{\sigma}_r) \, (1 - \hat{\sigma}_{r-1} \, \hat{\sigma}_r) \dots
(1 + (-1)^{r-k} \, \hat{\sigma}_k \, \hat{\sigma}_{k+1} \dots \hat{\sigma}_r) \,
f^{(\hat{\sigma})}_{k \to r+1} \; .
\ee
This is nothing but
the quantum version of the Zagier's factorization identities \cite{Zag}
(which can be obtained
from (\ref{zag}) by the projection $\hat{\sigma}_i \rightarrow q \, s_i$,
where $s_i$ are generators of the symmetric
group, $s_i^2 =1$). The proof of
identities (\ref{zag}) is straightforward (see e.g. \cite{isog2}).

\vspace{0.3cm}
We have the following result \cite{isog2}:

\vspace{0.2cm}
\noindent
{\bf Proposition 2.} {\it The explicit solution of the recurrent
relations (\ref{recur})
is given by the formula
\be
\label{solu}
\begin{array}{l}
(-1)^{r+1} \, \left( A X  A \right)^{i_1\dots i_r}_{j_1 \dots j_{r+1}}   \\ \\
\!\!\!\!\!\!\!\!
= \left( (1-R^2_{\underline{r}}) \left( 1 + R_{\underline{r-1}}
R^2_{\underline{r}} \right)
\cdots
\left( 1 + (-)^{r}  R_{\underline{1}} \dots R_{\underline{r-1}}
R_{\underline{r}}^2 \right)
\right)^{k_1\dots k_r \; 0}_{j_1 \dots j_r j_{r+1}}
\left( A_{1 \rightarrow r} \right)^{i_1\dots i_r}_{k_1\dots k_r} \; .
\end{array}
\ee
Here $\underline{1}, \dots , \underline{r-1}, \underline{r}$
label the copies of the vector space $V_{N+1}$;
$R_{\underline{r}} := R_{\underline{r,r+1}}$ and $\underline{r}$,
$\underline{r+1}$ are
the numbers of vector spaces $V_{N+1}$ where the $R$-matrix (\ref{rmat}) acts
nontrivially;
$i_1, \dots, j_1 \dots , k_1, \dots = 1,2, \dots N$ are the vector indices.
}

\vspace{.1cm}
\noindent
{\bf Proof}.
Using the definition of the antisymmetrizer (\ref{anti}) and the
braid relation (\ref{rell1})
we rewrite the right-hand side of (\ref{recur}) in the form
$$
\begin{array}{l}
A_{1 \rightarrow r+1} \, \left((-1)^{r} \,
\sigma_{r} \sigma_{r-1} \cdots  \sigma_{1} - {\bf 1} \right) \,
(XA)^{\langle  2\dots r |}_{|2 \dots r+1\rangle } \,  \\ \\
=  \left((-1)^{r} \, f_{1\to r+1} \,
\sigma_{r} \sigma_{r-1} \cdots  \sigma_{1} - \overline{f}_{1\to r+1} \right) \,
A_{2 \rightarrow r+1} \, (XA)^{\langle  2\dots r |}_{|2 \dots r+1\rangle } \,
\end{array} \ .
$$
Thus, the expression (\ref{solu}) is a solution of the equation (\ref{recur})
if the following identity holds
\be
\label{rel11}
\left(Y_{\underline{1 \to r+1}}
\right)^{\langle 1\dots r \; 0}_{1 \dots r+1\rangle } \, \overline{f}_{1 \to r} =
\left((-1)^{r+1} \, f_{1\to r+1} \,
 \sigma_{r} \cdots  \sigma_{1} + \overline{f}_{1\to r+1} \right)
\left( Y_{\underline{2 \to r+1}}
\right)^{\langle 2 \dots r \; 0}_{2 \dots r+1\rangle } \, ,
\ee
where
$$
Y_{\underline{k \to r+1}} =
(1-R^2_{\underline{r}}) \left( 1 + R_{\underline{r-1}} \,
R^2_{\underline{r}}
\right)
\cdots \left( 1 + (-)^{r} \, R_{\underline{k}} \,
R_{\underline{k+1}} \dots R_{\underline{r-1}} \,
R_{\underline{r}}^2 \right) \; .
$$

As indicated by the structure of indices, the vector space $V_{N+1}$ is a
direct sum of $V_N$ and a one-dimensional vector space $V_0$. Let
$P:V_{N+1}\rightarrow V_N$ and $1-P:V_{N+1}\rightarrow V_0$ be the
corresponding projectors. Similarly, let $\bar{P}$ and $1-\bar{P}$ be
the projectors associated to the dual space $\bar{V}_{N+1}$.

Define
$$
P_{1 \to r} :=  P_1\dots P_r\ ,\ P^{1\to r\; 0}=
\bar{P}_1\dots\bar{P}_r (1-\bar{P})_{r+1}\ .
$$
These operators split out the components which we need.

Eqn. (\ref{rel11}) can be rewritten now in the form
\be
\label{rel12}
\begin{array}{l}
P_{1 \to r+1} \, Y_{\underline{1 \to r+1}}
 \, \overline{f}^{(R)}_{\underline{1 \to r}} \,
P^{1\to r \; 0} \\ \\
= P_{1 \to r+1 } \, \left((-1)^{r+1} \, f^{(R)}_{\underline{1\to r+1}} \,
 R_{\underline{r}}  \cdots  R_{\underline{1}} +
 \overline{f}^{(R)}_{\underline{1\to r+1}} \right)
\, Y_{\underline{2 \to r+1}} \, P^{1 \to r \; 0}
\end{array}
\ee
(where $f^{(R)}_{\underline{1\to n}}$ and
$\overline{f}^{(R)}_{\underline{1\to n}}$
are the $(N+1)$-dimensional analogues of the operators (\ref{fbf}))
in view of relations
$$
\overline{f}^{(R)}_{\underline{1 \to r}} \, P^{1\to r}
= P^{1\to r} \, \overline{f}^{(\sigma)}_{1 \to r} \; ,
$$
$$
\begin{array}{lll}
Y_{\underline{2 \to r+1}} \, P^{1\to r \; 0}
&=& P^{1 \to r+1} \,
\left( P_{2 \to r+1 } \, Y_{\underline{2 \to r+1}} \, P^{2
\to r \; 0} \right)\\[1em]
&\equiv& P^{1 \to r+1} \, \left( Y_{\underline{2 \to r+1}}
\right)^{\langle 2 \dots r \; 0}_{2 \dots r+1\rangle }\ ,
\end{array}
$$
$$
P_{1 \to r+1 } \, \left((-1)^{r} \, f^{(R)}_{\underline{1\to r+1}} \,
R_{\underline{r}}  \cdots  R_{\underline{1}} -
\overline{f}^{(R)}_{\underline{1\to r+1}} \right) \, P^{1 \to r+1}
$$
$$
= \left((-1)^{r} \, f_{1\to r+1} \,
\sigma_{r} \sigma_{r-1} \cdots  \sigma_{1} -
\overline{f}_{1\to r+1} \right) \ ,
$$
which can be obtained directly using the explicit form of
the $R$-matrix (\ref{rmat1}).

The equation (\ref{rel12}) is fulfilled in view of the
identity (\ref{rel15}) which has been proven in Lemma 1. Indeed,
if in the identity (\ref{rel15}) we take the $R$-matrix representation
for the braid group $\rho_R(B_{M+1})$:
$\rho_R(\hat{\sigma}_k) = R_{\underline{k}}$
and act on this identity
from the left and right by the projectors
$P_{1 \to r+1 }$, $P^{1\to r \; 0}$, then
we deduce (\ref{rel12}). This
completes the proof of the Proposition 2. $\bullet$

\vspace{0.5cm}

In the next Section we need more information
about the $R$-matrix representation of the braid group just considered.
Introduce the Jucys-Murphy elements $J_r \in \rho_R(B_{M+1})$:
\be
\lb{jm}
J_1 = 1 \; , \;\;\; J_{r+1} = R_{\underline{r}} \, J_r \, R_{\underline{r}} \; .
\ee
These elements form a commuting set in $\rho_R(B_{M+1})$:
$[ J_r , \, J_m ] =0$
and satisfy
\be
\lb{be0}
R_{\underline{m}} \, J_r = J_r \, R_{\underline{m}} \; , \;\;\;   m < r \; .
\ee

\vspace{0.3cm}
\noindent
{\bf Proposition 3.}
{\it The explicit form of the components
$$Z_{r+1} := P_{1 \rightarrow r+1} (J_{r+1}) P^{1 \rightarrow r \; 0}$$
of Jucys-Murphy elements
(\ref{jm}) in terms of the structure constants $C_n$ and $\sigma_m$ is:
\be
\lb{be1}
\begin{array}{c}
Z_{r+1} = C_{r} + \sigma_r \, C_{r-1} \, \delta_{r+1}
+ \sigma_r \, \sigma_{r-1} \, C_{r-2} \, \delta_r \, \delta_{r+1} + \dots
\\[1em]
+ \sigma_r \, \sigma_{r-1} \dots \sigma_2 \, C_1 \, \delta_3 \dots
\delta_{r+1} \; .
\end{array}
\ee
The elements $Z_r$ have the following properties
\be
\lb{be3}
\sigma_{m} \, Z_{r} =  Z_r \, \sigma_m  \; , \;\;\; m \leq r-2 \; ,
\ee
\be
\lb{be2}
Z_{r+1} = C_r + \sigma_r \, Z_r \, \delta_{r+1} \; .
\ee
}

\vspace{0.3cm}
\noindent
{\bf Proof.}
Using the definition of Jucys-Murphy elements (\ref{jm}) and the representation (\ref{rmat1})
one can rewrite $P_{1 \rightarrow r+1} (J_{r+1}) P^{1 \rightarrow r \; 0}$
in the form
$$
\begin{array}{c}
P_{1 \rightarrow r+1} \, (R_{\underline{r}} \dots (R_{\underline{1}})^2 \dots R_{\underline{r}}) \,
 P^{1 \rightarrow r \; 0} \\[1em]
= P_{1 \rightarrow r+1} \,
\left( R_{\underline{r}} \, R_{\underline{r-1}} \, \dots R_{\underline{1}} \right)
 P^{0 \; 2 \rightarrow r+1} \, \delta_2 \delta_3 \dots \delta_{r+1}\ .
\end{array}
$$
Then, again using (\ref{rmat1}) we deduce (\ref{be1}).
To obtain the property (\ref{be3}) apply the operators $P_{1 \rightarrow r}$
and $P^{1 \rightarrow r-1 \; 0}$ to (\ref{be0}) from the left and
from the right (for $m \leq r-2$).
The recurrent relation (\ref{be2}) is a direct consequence
of the relations (\ref{be1}).  $\bullet$

\section{Relation to bar chain complex}

We recall first some standard material.
For every associative algebra ${\mathcal U}$ one forms the
bar complex $(C_n({\mathcal U}), b)$:
\begin{itemize}
\item[a)]  $C_n({\mathcal U})$ is the space
$\underbrace{{\mathcal U} \otimes \dots \otimes {\mathcal U} }_{n}$;
\item[b)] for
\be
T=a_1 \otimes a_2 \otimes \dots \otimes a_{n}
\in C_n({\mathcal U})\ ,
\lb{cht}\ee
the action of the boundary operator $b:
C_{n+1}({\mathcal U}) \rightarrow C_n({\mathcal U})$ is given by
\be
\lb{ho1}
\!\!\!\!\begin{array}{l}
b \, (a_1 \otimes a_2 \otimes \dots \otimes a_{n+1})  =
{\displaystyle \sum^n_{i=1}}
(-1)^{n-i} \,  (a_1 \otimes \dots \otimes a_i \, a_{i+1} \otimes \dots
\otimes a_{n+1} )  \; .
\end{array}
\ee
\end{itemize}
Using the associativity of ${\mathcal U}$ one can check directly that
$b^2 \, (a_1 \otimes a_2 \otimes \dots \otimes a_{n+1}) = 0$. Moreover,
for a unital algebra ${\mathcal U}$, the complex
$(C_n({\mathcal U}), b)$ is acyclic since one can define a homotopy operator
$\delta$: $C_{n}({\mathcal U}) \rightarrow C_{n+1}({\mathcal U})$
$$
\delta(a_1 \otimes a_2 \otimes \dots \otimes a_{n}) =
(a_1 \otimes a_2 \otimes \dots \otimes a_{n} \otimes 1)
$$
which satisfies $\delta \, b + b \, \delta = id$.

Let ${\mathcal U}$ be
the q-Lie algebra with
generators $\chi_i$ and defining relations (\ref{qla}).
Consider the chains $T$ (as in (\ref{cht})) with all elements
$a_k$ belonging to the set $\{ \chi_j \}$ of
the generators of the q-Lie algebra
and define their q-antisymmetric linear combination
which is given by the contraction with the numerical tensor
$A_{1 \rightarrow n} \in End(V^{\otimes n})$ (\ref{anti}):
\be
\chi_{1\rangle } \wedge \chi_{2\rangle } \wedge \dots \wedge \chi_{n\rangle }
:=A_{1 \rightarrow n} \,
\{ \chi_{1\rangle } \otimes \chi_{2\rangle } \otimes \dots \otimes \chi_{n\rangle } \}
\in C_{n}({\mathcal U})\ .
\ee
In fact we shall need the following chains:
\be
\lb{antc}
a \otimes \{ \chi_{1\rangle } \wedge \chi_{2\rangle } \wedge \dots \wedge \chi_{n\rangle } \}
\in C_{n+1}({\mathcal U})
\ee
with arbitrary $a\in {\mathcal U}$.
Denote the subspace of such chains by
$C_{n+1}({\mathcal U}, \wedge {\mathcal L}) \in$
${\mathcal U} \otimes \underbrace{{\mathcal L} \wedge \dots \wedge {\mathcal L}}_n$,
where ${\mathcal L}$ is a linear span of $\{ \chi_i \}$.

Now we formulate one of the main results of this paper.

\vspace{0.5cm}
\noindent
{\bf Proposition 4.}  {\it The boundary operator $b$ (\ref{ho1})
preserves the q-antisymmetry,
\be
\lb{hoho}
b: \;\; C_{n+1}({\mathcal U}, \wedge {\mathcal L}) \longrightarrow
C_{n}({\mathcal U}, \wedge {\mathcal L}) \; .
\ee
In other words, we have a complex
$(C_n({\mathcal U}, \wedge {\mathcal L}), b)$
which is a subcomplex of the bar com\-plex $(C_n({\mathcal U}), b)$.
}

\vspace{0.3cm}
\noindent
{\bf Proof.} Acting by the boundary operator (\ref{ho1}) on the
q-antisymmetric chains (\ref{antc}), we deduce
\be
\lb{ho3}
\begin{array}{l}
b \, \{ a \otimes \chi_{1\rangle } \wedge \chi_{2\rangle } \wedge \dots \wedge
\chi_{n+1\rangle } \}  =
(-1)^n \, A_{1 \rightarrow n+1} \,   \{ a \, \chi_{1\rangle } \otimes \chi_{2\rangle }
\otimes \dots \otimes \chi_{n+1\rangle } \}  \\ \\
+ a \otimes A_{1 \rightarrow n+1} \, \left[
{\displaystyle \sum^n_{k=1}} (-1)^{n-k} \,
\{ \chi_{1\rangle } \otimes \dots \otimes \chi_{k\rangle } \, \chi_{k+1\rangle } \otimes
\dots \otimes \chi_{n+1\rangle } \}  \right]
\end{array}
\ee
\be
\lb{hoch3}
\begin{array}{c}
= (-1)^n \, \overline{f}_{1 \to n+1} \,
\{ a \, \chi_{1\rangle } \otimes \chi_{2\rangle }  \wedge \dots
\wedge \chi_{n+1\rangle } \} \\ \\
+ a \otimes A_{1 \rightarrow n+1}  \,
{\displaystyle \sum^n_{k=1}} (-1)^{n-k} \, t^{\langle k}_{k, k+1\rangle }
\{  \chi_{1\rangle } \otimes \dots  \otimes \chi_{k\rangle } \otimes \chi_{k+2\rangle } \otimes
\dots \otimes \chi_{n+1\rangle } \} \; ,
\end{array}
\ee
where structure constants $t^i_{jk}$ satisfy
\be
\lb{ct}
C_{k k+1\rangle }^{\langle k} = (1-\sigma_k) \, t_{k k+1\rangle }^{\langle k} \; ,
\ee
tensors $\overline{f}_{1 \to n+1}$ are introduced in (\ref{anti}) -- (\ref{fbf})
and we have used the defining relations (\ref{qla}) and the properties of "antisymmetrizers"
$$
A_{1 \rightarrow n+1} \chi_{k\rangle } \, \chi_{k+1\rangle }=
A_{1 \rightarrow n+1} t_{k k+1\rangle }^{\langle k} \, \chi_{k\rangle } \; .
$$
Thus, to prove (\ref{hoho}) we need to represent the second term in (\ref{hoch3})
in the form
\be
\lb{hoch3a}
 a \otimes W_{1 \dots n+1\rangle }^{\langle 1 \dots n}  \,
 \{ \chi_{1\rangle } \wedge \dots  \wedge \chi_{n\rangle } \} =
 a \otimes W_{1 \dots n+1\rangle }^{\langle 1 \dots n}  \,
 A_{1 \to n} \, \{ \chi_{1\rangle } \otimes \dots  \otimes \chi_{n\rangle } \} \; .
\ee
Therefore it is sufficient to prove:

\vspace{0.3cm}
\noindent
{\bf Lemma 2.}
{\it There exist tensors $W_{n+1} := W_{1 \dots n+1\rangle }^{\langle 1 \dots n}$ such that
\be
\lb{x}
A_{1 \rightarrow n+1}  \, \left( \sum^n_{k=1} (-1)^{n-k} \, t^{\langle k}_{k, k+1\rangle } \,
\delta_{k+2 \rightarrow n+1} \right) =
W_{n+1}  \, A_{1 \rightarrow n}
\ee
where $\delta_{n+1 \to n} :=1$, $\delta_{k+2 \to n+1} :=
\delta_{k+2} \cdots \delta_{n+1}$ is a shift operator
$i \to i+1$ for $k < i \leq n$.
The explicit form of the tensors
$W_n$ (in terms of the quantum Lie algebra structure
constants $\sigma^{ij}_{kl}$ and $C^i_{jk}$) is:
$$
 W_{n+1} =  C_n -
(1-  \sigma_{n}) \, C_{n-1} \, \delta_{n+1} +
(1-\sigma_{n-1} +  \sigma_n \, \sigma_{n-1}) \,
C_{n-2} \, \delta_n  \, \delta_{n+1} - \dots
$$
\be
\lb{YZ}
\dots + (-1)^{n+1} \, [ 1-\sigma_{2} +  \sigma_3 \, \sigma_{2} - \dots + (-1)^{n+1}
\sigma_n \cdots \sigma_{2} ] \,
C_{1} \, \delta_3  \cdots \delta_{n+1}
\ee
$$
= Z_{n+1} - Z_n \, \delta_{n+1} +  Z_{n-1}  \, \delta_{n} \, \delta_{n+1} + \dots +
(-1)^{n+1} \, Z_{2}  \, \delta_{3} \dots \delta_{n+1}\ ,
$$
where the elements $Z_k$ are defined in (\ref{be1}), (\ref{be2}) (we fix
$Z_2 = C_1 =  W_2$).
}

\vspace{0.5cm}
\noindent
{\bf Proof.}
To prove (\ref{x}) we use the identities (\ref{rell1}) -- (\ref{rell4})
for the structure constants and the relation (\ref{ct})
between $C^j_{ik}$ and $t^j_{ik}$.

Eq. (\ref{x}) is obviously fulfilled for $n=1$.
Assume that (\ref{x}) is correct for $n=m-1$. We have to prove
(\ref{x}) for $n = m$. For the left hand side of (\ref{x}) we obtain:
$$
A_{1 \rightarrow m+1}  \, \left( \sum^{m}_{k=1} (-1)^{m-k} \, t^{\langle k}_{k, k+1\rangle } \,
\delta_{k+2} \cdots \delta_{m+1} \right)
$$
$$
 = A_{1 \rightarrow m+1}  \,
\left( \left[\sum^{m-1}_{k=1} (-1)^{m-k} \, t^{\langle k}_{k, k+1\rangle } \,
\delta_{k+2} \dots  \delta_m \right] \delta_{m+1}
+  t^{\langle m}_{m, m+1\rangle } \right)
$$
$$
 = \left(  -  f_{m+1} \, W_{m} \,  \delta_{m+1}
+  f_{m+1} \, f_m  \,
t^{\langle m}_{m, m+1\rangle } \right)  \, A_{1 \rightarrow m-1}
$$
$$
= \left[  - f_{m+1} \, W_{m}  \,  \delta_{m+1}
 +   x_{m+1} \,
C_{m} \right]  \, A_{1 \rightarrow m-1}\ .
$$
In the second equality we used (\ref{anti}), (\ref{x}); in the third
equality we used (\ref{ffxa}) with $n-m=2$.
Here and below, to simplify notation, we put
$f_m := f_{1 \to m}$ and $x_m := x^{(2)}_m$.

Comparing with the right-hand side of (\ref{x}), we obtain the equation
\be
\lb{equa}
- f_{m+1} \, W_{m}  \,  \delta_{m+1}
 +   x_{m+1} \, C_{m} = W_{m+1}  \, f_m \; ,
\ee
which is an identity for $m=2$ in view of (\ref{rell3}) and (\ref{rell4}).
We prove (\ref{equa}) by induction.
For this we note (formula (\ref{YZ})) that one has the inductive relation
\be
\lb{induy}
W_{n+1} = Z_{n+1} - W_n \, \delta_{n+1} \; .
\ee
Substitute to the left-hand side of
(\ref{equa}) the recurrence relations (\ref{fk}), (\ref{indux}), (\ref{induy})
and take into account  relations (\ref{rell3}), (\ref{be2}), (\ref{be3})
and the base of induction:
$$
- W_{m}  \,  \delta_{m+1}
+ f_m \, \sigma_m \, (Z_m - W_{m-1} \, \delta_m) \, \delta_{m+1} +
 (f_m + x_m \, \sigma_m \sigma_{m-1}) \, C_m
$$
$$
= - W_{m}  \,  \delta_{m+1} +
f_m \, ( \sigma_m \, Z_m \, \delta_{m+1} +  C_m )
+ W_m \, f_{m-1} \, \delta_{m+1}  \, \sigma_{m-1}
$$
$$
= ( Z_{m+1} - W_m \, \delta_{m+1}) -  f_{m-1} \, \sigma_{m-1} \, Z_{m+1}
+ W_m \, \delta_{m+1} \, f_{m-1} \, \sigma_{m-1}
$$
$$
= ( Z_{m+1} - W_m \, \delta_{m+1})(1 - f_{m-1} \, \sigma_{m-1})
= W_{m+1}  \, f_m \; .
$$
Thus, we proved (\ref{equa}) and, therefore, (\ref{x}). $\bullet$ \\
This finishes the proof of the Proposition 4. $\bullet$

\vspace{0.5cm}

The identity (\ref{x}) guarantees that the q-antisymmetric
chains (\ref{antc}) form the subcomplex in the bar complex for q-Lie algebra
${\mathcal U}$ (\ref{qla}). The identity
$$
(b)^2 \, \{ a \otimes \chi_{1\rangle }
\wedge \chi_{2\rangle } \wedge \dots \wedge \chi_{n+1\rangle } \} =0
$$
for the boundary operator (\ref{ho3}) directly follows from the
identity $b^2 = 0$ for the bar differential $b$ (\ref{ho1}).

In analogy with the theory of usual Lie algebras
(see e. g. \cite{CE}), we relate the com\-plex
$(C_{n}({\mathcal U}, \wedge {\mathcal L}),b)$
with the standard complex for the quantum Lie algebra.
Namely,  we relate the chains
$$
C_n = a \otimes A_{1 \rightarrow n} (\chi_{1\rangle } \otimes \dots \otimes \chi_{n})
\in C_{n}({\mathcal U}, \wedge {\mathcal L})
$$
with wedge polynomials in $\gamma$'s with coefficients in ${\mathcal U}$
\be
\label{mong}
M_n = a \, \gamma_{1\rangle } \wedge \gamma_{2\rangle }
\wedge \dots \wedge \gamma_{n\rangle } \; , \;\;\;
a \in {\mathcal U} \; ,
\ee
by means of the one-to-one map ${\bf i}: \; C_n \leftrightarrow M_n$,
\be
\lb{antc2}
{\bf i}(a \, \{ \gamma_{1\rangle } \wedge \gamma_{2\rangle }
\wedge \dots \wedge \gamma_{n\rangle } \}) =
 a \otimes  A_{1 \rightarrow n} \,
\{ \chi_{1\rangle } \otimes \chi_{2\rangle } \otimes \dots \otimes \chi_{n\rangle } \}
\ee
The BRST operator $Q$ (\ref{brst}), (\ref{qrr}), (\ref{solu})
acts on the expressions (\ref{mong}) from the right and,
according to (\ref{go}),
the variables $\Omega$'s in $Q$ should be considered as
"right $q$-derivatives" over $\gamma$'s
(in particular $\forall a \in {\cal U}$
the right action of $\Omega^i$ on $a$ gives zero).

\vspace{0.3cm}
\noindent
{\bf Proposition 5.}
\begin{itemize}
\item[1.] $Q^2 =0$;

\item[2.] {\it The map ${\bf i}$ (\ref{antc2}) is the chain map,
$b \circ {\bf i} = {\bf i} \circ Q$ }.
\end{itemize}

\vspace{0.3cm}
\noindent
{\bf Proof.} Consider
only first two terms in the expression (\ref{brst}):
 \be
 \lb{brst2}
 Q = \Omega^{\langle 1} \, \chi_{1\rangle } -
 \Omega^{\langle 2} \wedge \Omega^{\langle 1} \, t^{\langle 1}_{12\rangle } \, \gamma_{1\rangle }
 + \dots  \; .
  \ee
  Using (\ref{go}), (\ref{cross}),
(\ref{be1}) and (\ref{be2})
one can prove by induction the following relations:
$$
 \gamma_{1\rangle} \wedge \dots \wedge \gamma_{n\rangle} \, \Omega^{\langle n} =
 (-1)^n \,
  \Omega^{\langle 0} \, \sigma^{-1}_{0} \dots \sigma^{-1}_{n-1} \,
  \gamma_{0\rangle} \wedge \dots \wedge \gamma_{n-1\rangle}
$$
 \be
   \lb{app}
 + (-1)^{n-1} \,  \overline{f}^{(\sigma)}_{1 \to n} \,
 \sigma^{-1}_{1} \dots  \sigma^{-1}_{n -1}  \, \gamma_{1\rangle} \wedge \dots \wedge
 \gamma_{n-1\rangle} \, ,
  \ee
$$
\gamma_{1 \rangle} \wedge \cdots \wedge \gamma_{n-1 \rangle} \, \chi_{n \rangle} =
\sigma_{n-1} \dots \sigma_1 \, \chi_{1 \rangle}
\gamma_{2 \rangle} \wedge \cdots \wedge \gamma_{n \rangle} + Z_n \,
\gamma_{1 \rangle} \wedge \cdots \wedge \gamma_{n-1 \rangle} \; .
$$
 Thus, using (\ref{ffxb}) and (\ref{app}), the right action of
the first two terms (\ref{brst2}) of $Q$
 on (\ref{mong}) gives
\begin{equation}
\label{q2}
\begin{array}{c}
Q \left( a \,  \gamma_{1\rangle} \wedge \dots \wedge \gamma_{n\rangle} \right) =
(-1)^{n-1} \,  \overline{f}^{(\sigma)}_{1 \to n} \,
\left( a   \,  \chi_{1 \rangle} \gamma_{2\rangle} \wedge \dots \wedge
 \gamma_{n \rangle} \right) +
\\ \\
+ (-1)^{n-1}  \overline{f}^{(\sigma)}_{1 \to n} \,
  \sigma^{-1}_{1} \dots  \sigma^{-1}_{n -1}  \, Z_n \,
 \left(  a   \, \gamma_{1\rangle} \wedge \dots \wedge
 \gamma_{n-1\rangle} \right) -
\\ \\
- x^{(n-2)}_{ n} \,
 ( \sigma^{-1}_{1} \dots  \sigma^{-1}_{n -1}  ) \,
 ( \sigma^{-1}_{1} \dots  \sigma^{-1}_{n -2} ) \, C_{n-1} \,
 \left(  a   \, \gamma_{1\rangle} \wedge \dots \wedge
 \gamma_{n-1\rangle} \right) + \dots \; .
 \end{array}
\end{equation}
According to the Proposition 1 the first two terms in (\ref{brst2}),
together with the condition
\begin{equation}
\label{q2a}
Q^2|_{\rm terms \; linear \; in \; \chi}=0 \; ,
\end{equation}
define BRST operator $Q$ uniquely.
So we need only to check that the formula (\ref{q2}) (for $n=1,2$) is compatible with
$b \circ {\bf i} = {\bf i} \circ Q$ and eqs. (\ref{hoch3}), (\ref{hoch3a}), (\ref{antc2}).
It can be done directly. Therefore $Q^2 =0$ and
${\bf i}$ is a chain map. $\bullet$

The fact
that the map ${\bf i}$ (\ref{antc2}) is the chain map demonstrates the relation
between the standard complex for the q-Lie algebra
(with our BRST operator) and the subcomplex
$(C_{n}({\mathcal U}, \wedge {\mathcal L}),b)$ of the bar
complex $(C_{n}({\mathcal U}),b)$.


\vspace{0.3cm}
\noindent
{\bf Remark.}
One can show that
$$
W_{n+1} = P_{ 1 \to n+1} \, (\overline{f}^{(R)}_{\underline{1 \to n+1}} -1) \,
R_{\underline{1}} \dots R_{\underline{n}} \, P^{1 \to n,0} \; .
$$
The relation (\ref{equa}) can be obtained from the identity
\be
\lb{rfbr}
f^{(R)}_{\underline{1 \to m+1}} \,
(\overline{f}^{(R)}_{\underline{1 \to m}} \, R_{\underline{1}} \dots R_{\underline{m-1}})
R_{\underline{m}} =
(\overline{f}^{(R)}_{\underline{1 \to m+1}} \, R_{\underline{1}} \dots R_{\underline{m}} ) \,
f^{(R)}_{\underline{1 \to m}}
\ee
by action of $P^{1 \to m,0}$ and $P_{1 \to m+1}$
from the right and left. The identity (\ref{rfbr})
is equivalent to
$$
f^{(R)}_{\underline{1 \to m+1}} \,
\overline{f}^{(R)}_{\underline{1 \to m}}  =
\overline{f}^{(R)}_{\underline{1 \to m+1}} \,
f^{(R)}_{\underline{2 \to m+1}}
$$
which is nothing but the associativity condition
$x^{(1)}_{m+1} \, x^{(m-1)}_{m} = x^{(m)}_{m+1} T_1 x^{(1)}_{m} T_1^{-1}$
considered in the Remark in Section 3.
To prove all these statements we need the following
relations
$$
\begin{array}{l}
(\overline{f}^{(R)}_{\underline{1 \to n+1}} -1) \,
R_{\underline{1}} \dots R_{\underline{n}} \, P^{1 \to n,0}  \\ \\
= P^{ 1 \to n+1} \, W_{n+1} + \sum_{k=1}^n \, (-1)^k \, P^{1 \to k} (1 - P^{k+1}) \,
P^{k+2 \to n+1} \delta_{k+2} \dots \delta_{n+1}
\end{array}
$$
and
$$
f^{(R)}_{\underline{1 \to m}} \, P^{1 \to m} =
P^{1 \to m} \, f^{(\sigma)}_{1 \to m}\ .
$$

\section{A generalization}

{\bf 6.1.} We consider the similarity transformation of the $R$-matrix
\be
\lb{symtr}
R_{\underline{12}} \rightarrow
R'_{\underline{12}} = {\mathcal D}_{\underline{1}} \,
{\mathcal D}_{\underline{2}} \, R_{\underline{12}} \, {\mathcal D}^{-1}_{\underline{1}} \,
{\mathcal D}^{-1}_{\underline{2}}
\ee
which corresponds to the linear transformation of the generators
$L_B \rightarrow {\mathcal D}^C_B \, L_C$. The new matrix (\ref{symtr})
satisfies Yang-Baxter equation (\ref{rrr}) for all numerical matrices
${\mathcal D}$. If the matrix ${\mathcal D}$ has the special form
$$
{\mathcal D}^0_0 = 1 \; , \;\;\; {\mathcal D}^0_i =
{\mathcal D}^i_0 = 0 \; , \;\;\; {\mathcal D}^i_j = D^i_j \; ,
$$
then the new $R$ matrix $(\ref{symtr})$ has the same form (\ref{rmat})
with transformed structure constants
\be
\lb{strtr}
  \sigma_{12}  \rightarrow   \sigma'_{12} =
  D_1 \, D_2 \,  \sigma_{12} \, D^{-1}_1 \, D^{-1}_2 \; , \;\;\;
C^{\langle 2|}_{|12\rangle} \rightarrow {C \, '}^{\langle 2|}_{|12\rangle} =
D_1 \, D_2 \, C^{\langle 2|}_{|12\rangle}  \, D^{-1}_2 \; .
\ee
This transformation leaves invariant the linear space (denoted by
${\mathcal L}$) spanned by the $N$ elements $\chi_i$.
The q-Lie algebra is quadratic for any
choice of the basis in ${\mathcal L}$ but the structure constants in
(\ref{qla}) change according to (\ref{strtr}).

The interesting special case is when the
matrix $R$ is conserved under transformation
(\ref{symtr}) $R_{\underline{12}} = R'_{\underline{12}}$.
Then the $N \times N$ matrix $D^i_j$ is such that
\be
\lb{rsym}
D_1 \, D_2 \,  \sigma_{12}  =  \sigma_{12} \, D_1 \, D_2 \; , \;\;\;
 D_1 \, D_2 \, C^{\langle 2|}_{|12\rangle} = C^{\langle 2|}_{|12\rangle} \, D_2 \; .
\ee

\vskip .3cm
\noindent{\bf Lemma.} {\it In this case the $R$-matrix
${\mathcal D}_{\underline{1}} \, R_{\underline{12}} \,
{\mathcal D}^{-1}_{\underline{1}}$ satisfies Yang-Baxter equation.}

\vskip .3cm
The proof is straightforward.
Now one can consider a quadratic algebra with defining relations
$$
\left( {\mathcal D}_{\underline{1}} \, R_{\underline{12}} \,
{\mathcal D}^{-1}_{\underline{1}} \right)
\, L'_{\underline{1} \rangle}  \, L'_{\underline{2} \rangle} =
L'_{\underline{1} \rangle}  \, L'_{\underline{2} \rangle} \; .
$$
In the components, for the quantum vector $L'_A = \{ \xi_0, \xi_i \}$, we obtain
\be
\lb{qla2}
\xi_i  \, \xi_0 =  D^j_i \, \xi_0 \, \xi_j  \ \ {\rm and} \ \
(1 - D_1 \, \sigma_{12} \, D^{-1}_1) \, \xi_{1\rangle } \, \xi_{2\rangle } =
D_1 \, C_{12\rangle }^{\langle 1} \, \xi_0 \, \xi_{1\rangle }\ .
\ee
This algebra is related to the q-Lie algebra (\ref{qla}) by a transformation
\be
\lb{qla3}
\chi_i = \xi_0^{-1} \xi_i\ .
\ee

\vspace{0.5cm}
{\bf 6.2.} There is a generalization of the exterior algebra
(\ref{owedge}) --(\ref{go}) and (\ref{cross})
arising from the considerations just above.
One can introduce a grading operator $g$ with commutation relations
\be
\lb{grg}
g \, \gamma_{1\rangle} = D_1 \, \gamma_{1\rangle} \, g \; , \;\;\;
g \, \Omega^{\langle 1} = \Omega^{\langle 1} \, D^{-1}_1 \, g \; , \;\;\;
g \, \chi_{1\rangle} = D_1 \, \chi_{1\rangle} \, g \; ,
\ee
where $D^i_j$ is the matrix which satisfies (\ref{rsym}).
Note that the grading operator $g$ can be included into play if we
consider the quadratic algebra (\ref{qla2}) and take (\ref{qla3}) and
$g = \xi_0^{-1}$.

Relations (\ref{grg}) are consistent
with relations (\ref{qla}),
(\ref{owedge}) --(\ref{go}) and (\ref{cross}).
Then the Woronowicz differential algebra
(\ref{owedge}) --(\ref{go}), (\ref{cross})
can be rewritten for the new generators
$$
\tilde{\Omega}^{\langle 1} = \Omega^{\langle 1} \, g^{-1} \; , \;\;\;
\tilde{\gamma}_{1\rangle} = g \, \gamma_{1\rangle}
$$
in the form:
\be
\lb{go2}
 \tilde{\gamma}_{2\rangle} \, \tilde{\Omega}^{\langle 2} = - \tilde{\Omega}^{\langle 1} \,
D^{-1}_1 \, \sigma^{-1}_{12} \, D_1 \, \tilde{\gamma}_{1\rangle} + I_2  \; ,
\ee
 \be
 \lb{owedge2}
\tilde{\Omega}^{\langle 1} \wedge \tilde{\Omega}^{\langle 2} \wedge
\dots \wedge \tilde{\Omega}^{\langle n}
= \tilde{\Omega}^{\langle 1} \otimes \tilde{\Omega}^{\langle 2} \otimes \dots \otimes
\tilde{\Omega}^{\langle n} \, A^{(D)}_{1 \rightarrow n} \; ,
 \ee
  \be
 \lb{gwedge2}
\tilde{\gamma}_{1 \rangle} \wedge \tilde{\gamma}_{2\rangle} \dots
\wedge \tilde{\gamma}_{n\rangle} =
A^{(D)}_{1 \rightarrow n}
\, \tilde{\gamma}_{1\rangle} \otimes \tilde{\gamma}_{2\rangle} \dots
\otimes \tilde{\gamma}_{n\rangle} \; ,
 \ee
  \be
\lb{cross2}
 \chi_{2 \rangle} \, \tilde{\Omega}^{\langle 2} =  \tilde{\Omega}^{\langle 1} \, \left(
 \sigma_{12} D_1 \, \chi_{1 \rangle} +  C^{\langle 2}_{ 12\rangle} \right) \, ,
 \;\; \tilde{\gamma}_{1 \rangle} \, \chi_{2\rangle} = \sigma_{12} D_1 \, \chi_{1\rangle} \,
\tilde{\gamma}_{2\rangle} + C^{\langle 1}_{12\rangle} \, \tilde{\gamma}_{|1\rangle} \; ,
 \ee
where
$A^{(D)}_{1 \rightarrow n} =
D_1^{-1} \cdots D_{n-1}^{-1} \, A_{1 \rightarrow n} \, D_1 \cdots D_{n-1}$.

This generalized algebra embeds into the bar complex; the embedding
is similar to the one in the particular case $D=1$ and we don't give details.

\vspace{0.5cm}
{\bf 6.3.} The quantum Lie algebras defined by eqs (\ref{qla}),
(\ref{rell1})-(\ref{rell4}) generalize the
usual Lie (super)algebras. Indeed, in the non-deformed case,
the braid matrix $\sigma^{mk}_{ij} = (-1)^{(m)(k)} \, \delta^m_j \, \delta^k_i$
is a super-permutation matrix (here $(m)=0,1$ is the parity of a generator
$\chi_m$).
Eqn. (\ref{rell1}) is fulfilled (and we have additionally
$\sigma^2 =1$). Eqn.
(\ref{rell2}) coincides with the Jacobi identity for Lie (super)-algebras.
Eqn. (\ref{rell3}) is then equivalent to the
$Z_2$-homogeneity condition $C^i_{jk} =0$ for $(i) \neq (j)+(k)$.
Eqn. (\ref{rell4}) follows from (\ref{rell3}) in this case.

In this sense the exterior algebras (\ref{owedge}) --(\ref{go}),
(\ref{go2}) -- (\ref{gwedge2})
generalize the Heisenberg algebras of the fermionic (and bosonic) ghosts and
anti-ghosts. The grading matrix $D$ now is the matrix
$D^m_n = (-1)^{(n)} \delta^m_n$ which defines
the well known automorphism of superalgebras, and we have
$$
(D_1\sigma_{12}D_1^{-1})^{mk}_{ij} = -(-1)^{((m)+1)((k)+1)}
\, \delta^m_j \, \delta^k_i \; .
$$

\section{Conclusion}

We have investigated
the special class of quadratic algebras ${\mathcal U}$, the
quantum Lie algebras (\ref{qla}). We have considered the
exterior extension of ${\mathcal U}$ by the ghost algebra
(\ref{owedge}), (\ref{gwedge}),
(\ref{go}), (\ref{cross}) and constructed the BRST operator $Q$.
Using the BRST operator (\ref{brst}), (\ref{qrr}), (\ref{solu}) one can
build analogues of the standard and de Rham complexes for quantum
Lie algebras.
As we have shown in this paper, there is a map ${\bf i}$ from
the standard complex to the subcomplex in the bar complex of
${\mathcal U}$. The map ${\bf i}$ induces (from the bar differential $b$) a differential
on the standard complex. We compared this induced differential with the
operator $Q$.
The operator $Q$ is defined uniquely by two properties:
the initial terms (\ref{brst2}) and the condition (\ref{q2a}).
We verified that the induced  differential satisfies these two properties.
Therefore it coincides with $Q$. Moreover it follows that $Q^2=0$.
For the construction of the map ${\bf i}$ we need to
check a number of complicated identities on the structure constants
$\sigma^{kl}_{ij}$ and $C_{ij}^k$. The similar identities have appeared
in the inductive definition of the BRST operator for q-Lie algebras.
An elegant proof of some of these identities is based on
combining the constants $\sigma^{kl}_{ij}$ and $C_{ij}^k$ into a
bigger matrix $R^{AB}_{CD}$ which realizes an $R$-matrix representation of
the braid group ${\mathcal B}_\bullet$.
Then we use the properties of the Jucys - Murphy elements in
${\mathcal B}_\bullet$, the quantum shuffle product and several
important identities in the braid group algebra.

We presented a generalization of the constructions above to the situation
when a q-Lie algebra is equipped with a grading operator.

Some results from Section 5 can be considered as an explicit realization of
certain facts from the paper \cite{BG} for the special choice of
the inhomogeneous quadratic algebras.

\vspace{0.3cm}

Note that the particular example of our construction of the BRST operator
for the case of the quantum Lie algebra $U_q(gl(N))$ has been already
considered in detail in our paper \cite{isog1}.
In this case we have constructed also the quantum analogues of
the anti-BRST operator and the quantum Laplace operator.

Also,we hope that our construction of the BRST operator for the quantum Lie
algebras will be
useful for constructions of BRST operators for any quadratic algebras
(even for infinite dimensional algebras, such as $W_3$)
and for producing the quantum $W$-algebras
with the help of the quantum affine algebras
(see \cite{Sev} and references therein) by means of
the quantum analogue of the Hamiltonian
reduction procedure via the quantum BRST technique.

\vspace*{-5pt}

\section*{Acknowledgements}
The authors are especially grateful to A. Gerasimov,
P. Etingof and A. Polischuk for useful discussions.
V. G. and A. P. I. also
thank the Max-Planck-Institut f\"{u}r Mathematik
in Bonn and Marseille University,
where the considerable part of this work was done, for their kind hospitality
and support. We also thank the referee for useful remarks.

The work of A.P.I. was supported in part by the RFBR (Grant No. 03-01-00781)
and INTAS (Grant No. 03-51-3350).

\vspace*{-5pt}


\end{document}